\documentclass[12pt, a4paper]{article}

\usepackage{euscript} 
\usepackage{a4wide} 
\usepackage[utf8]{inputenc} 
\usepackage[T2A]{fontenc} 
\usepackage{amsfonts} 
\usepackage{amssymb, amsthm} 
\usepackage{amsmath}
\usepackage{tikz}
\usepackage{tikz-cd}
\usepackage{graphicx} % Required for including images
\usepackage{floatrow}
\usepackage{caption}
\usepackage[all,cmtip]{xy}

\newtheorem*{theorem*}{Theorem}

\date{January 20, 2022}
\DeclareMathOperator{\Conj}{Conj}
\DeclareMathOperator{\Adconj}{Adconj}
\DeclareMathOperator{\FR}{FR}
\DeclareMathOperator{\FQ}{FQ}

\title{\textsc{On Schreier varieties of racks}}
\author{Georgii Kadantsev and Aleksandra Shutova} 

\usepackage[
backend=biber,
style=ieee,
url=false,
doi=false,
]{biblatex}

\begin{document}
\maketitle 
\begin{abstract} 
We prove that a subrack of a free rack is free and suggest a method to prove a similar statement about involutory racks.
\end{abstract}

Key words: free rack, quandle, Schreier varieties, Nilsen-Schreier theorem.

\section*{Introduction}
A rack is a set equipped with a two binary operations $(R, \triangleright, \triangleright^{-1})$ such that the following equalities hold for every $x,y,z\in R$:
\begin{itemize}
    \item[R1] $(x\triangleright y) \triangleright z = (x\triangleright z) \triangleright (y \triangleright z)$, \\
    $(x\triangleright^{-1} y) \triangleright^{-1} z = (x\triangleright^{-1} z) \triangleright^{-1} (y \triangleright z)$;
    \item[R2] $(x \triangleright y) \triangleright^{-1} y = x = (x \triangleright^{-1} y) \triangleright y$.
\end{itemize}
R1 states simply that the map $x \to x \triangleright y$ is an endomorphism of $Q$ for every $y \in Q$.
R2 implies that every such map is an automorphism.
A rack does not need to be associative or to have an identity. 

A rack with $x \triangleright x = x$ is called a quandle.
A rack in which $\triangleright = \triangleright^{-1}$ is called involutory. 
Involutory quandles have been studied extensively under different names (symmetric sets, symmetric groupoids, see \cite{History}).

Any group $G$ provides an example of a quandle $\Conj G$ with $x\triangleright y = y^{-1} x y$. 
$\Conj$ can be considered as a functor from the category of quandles to the category of groups. 
There exists a left adjoined functor to $\Conj$ which we denote by $\Adconj$. 
$\Adconj Q$ is the universal group in which to represent the quandle $Q$ as a set closed under conjugation.  
We call it the associated group of the quandle.

One strong motivation for studying quandles and racks is provided by knot theory.
There is a natural construction of a quandle $Q(K)$ for any knot $K$ using its diagram. 
It is called {\it the knot quandle} or {\it the fundamental quandle} of the knot (see \cite{joyce_classifying_1982} for details).
This construction gives a full invariant of knots and other invariants can be derived from it (see \cite{Kamada}, \cite{Sam}).
For example, the fundamental group of a knot is obtained as the associated group of its fundamental quandle. 

A variety is a class of algebraic structures of the same type satisfying a set of identities \cite{birkhoff_1935}.
A variety of algebras in which subalgebras of free algebras are free is called a Schreier variety.
So, by the Nielsen-Schreier theorem the variety of groups is Schreier.
Schreier varieties of linear algebras have been studied in \cite{umirbaev_schreier_1994, MoreSchreier}.

V. Bardakov, M. Singh and M. Singh in \cite[Problem 6.12]{Bardakov_2019} raised the question about an analogue of Nielsen–Schreier theorem for quandles: is it true that any subquandle of a free quandle is free.
It was answered affirmatively in \cite{Kad}. 
In this work we generalise this result to racks and propose a way to extend it to involutory racks:

\begin{theorem*}
    Any subrack of a free rack is free.
\end{theorem*}

\section*{Algebraic representation of racks and quandles}

Free rack on $X$ is a rack that satisfies the universal property: given any function $\rho \colon X \to R$, where $R$ is an arbitrary rack, there exists a unique homomorphism $\overline{\rho} \colon \FR(X) \to R$, such that $\overline{\rho} \circ \varphi = \rho,$ making the following diagram commute (here $\varphi \colon X \to \FR(X)$ is an embedding of $X$ into $\FR(X)$):

$$\xymatrix{
X \ar[dr]^{\phi} \ar[d]^\rho & \\
R & \ar[l]^{\overline{\rho}} \FR(X)
}$$

We will use the following construction of a free rack on $X$ \cite{fenn_racks_1992}. 
On the set $X \times F(X)$, where $F(X)$ is a free group, generated by $X$, we will define $\triangleright$ as follows:
$$(a, u) \triangleright (b, v) = (a, uv^{-1}bv)$$
$$(a, u) \triangleright^{-1} (b, v) = (a, uv^{-1}b^{-1}v)$$
for all $a, b \in X, \ u, v \in F(X)$. 

A free quandle $\FQ(X)$ on $X$ is a union of conjugacy classes of elements of $X$ in $F(X)$
with the operation difined the following way:
$x \triangleright y = x ^ y = y^{-1}xy \ \forall x, y \in FQ(X)$ \cite{joyce_classifying_1982}.

A free involutory quandle is a union of conjugacy classes of elements of $X$ in $\langle X |\ x^2 = 1\ \forall x\in X\rangle$.  

\section*{Proof}
For convenience we will denote  $((r_0 \triangleright ^ {\epsilon_1} r_1) \triangleright ^ {\epsilon_2} \dots) \triangleright^{\epsilon_{n}} r_n = r_0 \triangleright^{\epsilon_1} r_1 \triangleright ^ {\epsilon_2} \dots \triangleright ^ {\epsilon_{n}} r_n$.
We will also write $r^n, n \in \mathbb{Z}$ instead of $r \triangleright ^ {\epsilon} r \triangleright ^ {\epsilon} \dots \triangleright ^ {\epsilon} r$, where $\epsilon = sign(n)$.

Consider $f\colon \FR(X) \to \FQ(X), \ f((x,w)) = w^{-1}xw = x^w$.
It is clear that $f$ is a rack homomorphism.
We will consider an arbitrary subrack $R \subset \FR(X)$ and show that it is free.
Since $f$ is a homomorphism, the image of $R$ is a subquandle $Q \subset \FQ(X)$, and thus is free.

Then a basis $S_Q$ exists, such that $Q = \langle S_Q \rangle$. 
Any element of $Q$ can be represented as $q_0 \triangleright ^ {\epsilon_1} q_1 \triangleright ^ {\epsilon_2} q_2 \dots \triangleright ^ {\epsilon_n} q_n$, where $q_i \in S_Q$. 

The preimage of $x^w$ in $\FQ(X)$ is the subrack $\{ (x, x^n w) \ | n \in \mathbb{Z} \}$, generated by any one of its elements. 
To prove this, assume that $(x, w_1)$ and $(y, w_2)$ are such that $f((x,w_1)) = x^{w_1} = f((y,w_2)) = y^{w_2}$.
Since $x$ and $y$ belong to the same conjugacy class in $FR(X)$, $x=y$.
Now $x^{w_1} = x^{w_2}$ implies $x = x^{w_1w_2^{-1}}$, which is possible only if $w_1w_2^{-1} = x^n$ for some $n\in \mathbb{Z}$.

Note that $f^{-1}(x^w) \subset R$.

From each preimage of $q_i = x^w \in S_Q$ choose $r_i = (x, w^\prime)$, where $w^\prime$ does not start with a power of $x$.  
Note that every $r_i$ is unique, otherwise $q_i$ are not independent from each other and do not form a basis of $Q$. 
We will denote this set of $r_i$ by $S_R$ and show that it generates $R$ freely. 

For $r = (x, w)\in R$ consider $f(r) = x^w \in Q$. 
Sinсe $r$ is contained in the preimage of $x^w$, 
it can be represented as 
$$(q_1 \triangleright ^ {\epsilon_1} q_2 \triangleright ^ {\epsilon_2} \dots \triangleright ^ {\epsilon_{n- 1}} q_m)^n, \quad q_i\in S_Q, \ q_1 \neq q_2.$$

Using the equality 
$(r \triangleright^{\epsilon} t)^k = r^k \triangleright^{\epsilon} t,$
which holds in every rack, we obtain 
$$r = q_1^n \triangleright ^ {\epsilon_1} q_2 \triangleright ^ {\epsilon_2} \dots \triangleright ^ {\epsilon_{n- 1}} q_m$$

Now let us show that this representation is unique.
Assume that $x_0, x_1, \ldots x_n$ and $y_0, y_1,\ldots y_m$ are such that 
\[k,l\neq 0, x_0 \neq x_1, y_0 \neq y_1, \quad  x_0^k \triangleright^{\epsilon_1} x_1 \triangleright^{\epsilon_2} x_2 \ldots \triangleright^{\epsilon_n} x_n = 
 y_0^l \triangleright^{\xi_1} y_1 \triangleright^{\xi_2} y_2 \ldots \triangleright^{\xi_m} y_m \]
Denote $f(x_i)$ and $f(y_j)$ by $\overline{x_i}$ and $\overline{y_j}$ respectively.
Applying $f$ to both sides gives us 
\[\overline{x_0} \triangleright^{\epsilon_1} \overline{x_1} \triangleright^{\epsilon_2} \overline{x_2} \ldots \triangleright^{\epsilon_n} \overline{x_n} = 
 \overline{y_0} \triangleright^{\xi_1} \overline{y_1} \triangleright^{\xi_2} \overline{y_2} \ldots \triangleright^{\xi_m} \overline{y_m} \]
Since this is an equation on basis elements in $Q$, we have $n=m$, $\epsilon_i = \xi_i$ and $\overline{x_i} = \overline{y_i}$ for every $i$.
The mapping $f$ is injective on elements of $S_R$, which means that $\overline{x_i} =\overline{y_i}$ implies $x_i = y_i$.
Now all $x_i$, where $i \geq 1$, can be cancelled out.
What is left is
\[ x_0^k = x_0^l.\]
In $FR(X)$ this is possible only if $k = l$.
This concludes the proof.

A similar proof can be carried out with a construction of free involutory racks as the proof above does not change with $\triangleright = \triangleright^{-1}$, given that the variety of involutory quandles is Schreier:
\begin{theorem*}
    Every subrack of a free involutory rack is a free involutory rack. 
\end{theorem*}

\printbibliography

{\scshape Laboratory of Continuous Mathematical Education, Saint Petersburg, Russia\\}
\textit{E-mail address:} \texttt{kadantsev.georg@yandex.ru}
\vspace{0.2cm}

{\scshape Laboratory of Continuous Mathematical Education, Saint Petersburg, Russia\\}
\textit{E-mail address:} \texttt{beloshapkoa14@gmail.com}
\end{document}